\documentclass[reqno,12pt]{amsart}
\usepackage{amsfonts,amssymb,amsbsy,amsmath,amsthm}

\newtheorem{lemma}{Lemma}[section]
\newtheorem{theorem}[lemma]{Theorem} 
\newtheorem{follow}[lemma]{Corollary}
\newtheorem{pr}[lemma]{Proposition}

\theoremstyle{definition}

\newtheorem{example}[lemma]{Example}

\newtheorem{remark}[lemma]{Remark}

\newcommand{\bel}{\begin{equation} \label}
\newcommand{\ee}{\end{equation}}

\newcommand{\rd}{{\mathbb R}^{2}}
\newcommand{\re}{{\mathbb R}}

\newcommand{\q}{\quad}
 \newcommand{\curl}{\operatorname{curl}}
  \newcommand{\ran}{\operatorname{Ran}}
    \newcommand{\sgn}{\operatorname{sgn}}
     \newcommand{\dive}{\operatorname{div}}
       \newcommand{\grad}{\operatorname{grad}}
      \newcommand{\slim}{\operatorname{s-lim}}
         \newcommand{\supp}{\operatorname{supp}}

{

\let\goth\mathfrak
\let\Bbb\mathbb
\let\cal\mathcal

\begin{document}

\title[Translationally Invariant Magnetic   Operators]{On   Spectral Properties of  Translationally Invariant Magnetic Schr\"odinger Operators}

  
  \author{ D.   Yafaev}
\address{ IRMAR, Universit\'{e} de Rennes I\\ Campus de
  Beaulieu, 35042 Rennes Cedex, FRANCE}
\email{yafaev@univ-rennes1.fr}

\begin{abstract}
We consider a class of translationally invariant magnetic fields such that the corresponding potential has a constant direction.
 Our goal is to study    basic spectral
properties of the   Schr\"odinger operator ${\bf H}$ with such a potential. In particular, we show that the spectrum of ${\bf H}$ is
   absolutely continuous and we find its location. Then we study  the
long-time behaviour of solutions $\exp(-i {\bf H} t)f$ of the time dependent Schr\"odinger equation.
      It turnes out that a quantum particle remains localized in the plane orthogonal to the direction of the potential. Its propagation in this direction is determined by group velocities. It is to a some extent similar to a evolution of a one-dimensional free particle but ``exits" to $+\infty$ and $-\infty$ might be essentially different.
\end{abstract}

 \maketitle

\section{Introduction}
\setcounter{equation}{0}

{\bf 1.1.} 
Traslationally invariant magnetic fields $B(x)=(b_{1}(x), b_{2}(x), b_{3}(x))$, $x=(x_{1}, x_{2}, x_{3})$, $\dive B(x)=0$, give important examples where a non-trivial information can be obtained about spectral properties of the corresponding  Schr\"odinger operators ${\bf H} $.
We suppose for definiteness that  $B(x)$ does not depend on the $x_{3}$-variable so that ${\bf H} $ commute with translations along the $x_{3}$-axis. 
There are two   essentially different (and in some sense extreme) classes of traslationally invariant magnetic fields. 
 
 The first class consists of fields $B(x)=(0,0,b_{3}(x_{1}, x_{2}))$ of constant direction. For such fields, the momentum $p$ of a classical particle in the $x_{3}$-direction is conserved, and in the Schr\"odinger equation the variable $x_{3}$ can be separated. Thus, we arrive  to a two-dimensional problem in the $(x_{1}, x_{2})$-plane. Furthermore, if  $b_{3}$ is a function of $r=(x_{1}^2+ x_{2}^2)^{1/2}$ only, then we get a set of problems on the half-line $r>0$ labelled by the magnetic quantum number $m$. 
 The most important example of this type is a constant magnetic field $b_{3}(r)=const$ (see \cite{LL}). Some class of functions $b_{3}(r)$ decaying as $r\to\infty$ was discussed in  \cite{MS} (see also  \cite{CFKS}) where new interesting effects were found. Another famous   case is $b_{3}(x_{1}, x_{2})=\delta(x_{1}, x_{2})$ ($\delta(\cdot)$ is the Dirac delta-function) studied in \cite{AB}. Scattering by an arbitrary short-range (decaying faster than $|x|^{-2-\varepsilon}$, $\varepsilon>0$, as $|x|\to \infty$) magnetic field $b_{3}(x_{1}, x_{2})$ turns out to be rather similar to this particular case (see \cite{Yaf}).

The second class consists of fields $B(x)=(b_{1}(x_{1}, x_{2}),b_{2}(x_{1}, x_{2}), 0)$ orthogonal to the $x_{3}$-axis. In this case the corresponding magnetic potential
$A(x)$, defined (up to gauge transformations) by the equation $\curl A(x)=B(x)$,  can be chosen as
 \bel{w3}
 A(x)=(0,0,-a(x_{1}, x_{2}))  
 \ee
 so that it has the constant direction.
In contrast to fields of the first class, now the variable $x_{3}$ cannot be separated
in the Schr\"odinger equation. 
 Nevertheless due to the invariance with respect to translations along the $x_{3}$-axis the  operator (we always suppose that the charge of a particle is equal to $1$)
  \bel{MH}
{\bf H}  = (i\nabla + A (x))^2
\ee
can be realized, after the Fourier transform in the variable $x_{3}$,  in the space $L^2(\re;L^2(\re^2))$ as the operator of multiplication by the operator-valued function
 \bel{u1ax}
H (p) = -\Delta + (a(x_{1}, x_{2})+p)^2: L^2(\re^2)\to L^2(\re^2),\q p\in \re. 
\ee
 Moreover, if  $a(x_{1}, x_{2})=a(r)$, then
the
subspaces   with fixed magnetic quantum number $m \in {\mathbb Z}$ are invariant subspaces of $H (p)$ so that the operator
$H (p)$ reduces to the orthogonal sum over $m \in {\mathbb Z}$
of the operators
\bel{u1a}
H_m(p) = -\frac{1}{r} \frac{d}{dr}r  \frac{d}{dr} +
\frac{m^2}{r^2} + (a(r)+p)^2 
\ee
acting  in the space ${\cal H}=L^2(\re_+;rdr)$. 
 In this case the field  is given by the equation 
 \bel{magnfi}
 B(x)=b(r)(-\sin \theta, \cos \theta,0)
 \ee
  where $b(r)=a'(r)$ and $\theta$ is the polar angle. 
 Thus, vectors $B(x)$   are tangent to circles centered at the origin.
 An important example of such type is a field created by a current along an infinite straight wire (coinciding with the $x_{3}$-axis). In this case   $b(r)=b_{0} r^{-1}$ so that  $a(r)=b_{0} \ln r$.
 The Schr\"odinger operator with such magnetic potential was studied in \cite{Y}.

  \medskip

{\bf 1.2.} 
In this article we consider magnetic fields  \eqref{magnfi} with a sufficiently arbitrary function $b(r)$.
 Our goal is to study    basic spectral
properties of the corresponding Schr\"odinger operator
${\bf H}  $
 such as the absolute continuity, location and
multiplicity of the spectrum, as well as   the
long-time behaviour of the unitary group $\exp{(-i{\bf H}t)}$.
 We emphasize that for magnetic fields considered here, the problem is genuinely three-dimensional, and actually the motion of a particle in the $x_{3}$-direction is of a particular interest.
 

     Using the cylindrical invariance of field \eqref{magnfi}, we can start either from translational or from rotational (around the $x_{3}$-axis) symmetries. 
The rotational   invariance implies that the operator ${\bf H}$ is the orthogonal sum of its restrictions ${\bf H}_m$ on the
subspaces of functions   with   magnetic quantum number $m \in {\mathbb Z}=\{0,\pm 1, \pm 2,\ldots\}$.
  It can be identified  with the operator (we keep the same notation for this operator)
\bel{Hm1}
{\bf H}_m = -\frac{1}{r} \frac{d}{dr}r  \frac{d}{dr} +
\frac{m^2}{r^2} + (i\frac{d}{dx_{3}} +a(r) )^2 
\ee
acting in the space ${\goth H}=L^2(\re_+\times \re ;rdr dx_{3})$.  In view of   the translation invariance, 
every operator  ${\bf H}_m$  
can be realized (again after the Fourier transform in the variable $x_{3}$) in the space $L^2(\re;L^2(\re_{+};rdr))$ as the operator of multiplication by the operator-valued function    $  H_m(p) $ defined by \eqref{u1a}. 

Suppose now that $b(r)$ does not tend to zero too fast so that  
    $a(r) \to \infty$ as $r \to \infty$. Then the spectrum of each operator $H_{m}(p)$ is discrete. Let
$  \lambda_{n,m}(p)$, $n \in {\mathbb N}=\{1,2,\ldots\}$, be the
increasing sequence of its eigenvalues (they are simple and positive), and let $ \psi_{n, m}(r,p)$ be the corresponding sequence of its eigenfunctions.
  The functions $\lambda_{n,m}(p)$ are known as dispersion curves of the problem. They determine
  the spectral properties of the operator ${\bf H}_{m}$.
  
   Note that if $a(r)$ is replaced by 
  $-a(r)$, then  $\lambda_{n,m}(p)$ is replaced by $-\lambda_{n,m}(p)$, so that the case
  $a(r) \to -\infty$ as $r \to \infty$ is automatically included in our considerations.
   Recall that, for a magnetic field $B(x)$,  the magnetic potential $A(x)$ such that $\curl A(x) =B(x)$ is defined up to a gauge term $\grad \varphi(x)$. In particular for magnetic fields \eqref{magnfi} in the class of potentials $A(x)=(0,0, - a(r))$ one can always add to $a(r)$ an arbitrary constant $c$. This leads to the transformations $\lambda_{n,m}(p)\mapsto
\lambda_{n,m}(p-c)$ and $\psi_{n,m}(r,p)\mapsto
\psi_{n,m}(r,p-c)$.

   \medskip

{\bf 1.3.} 
 The precise definitions of the operators ${\bf H}_{m}$ and ${\bf H} $ and  their decompositions into the  direct integrals over the operators $H_{m}(p)$ and $H (p)$  are given in   Section 2. To put it differently, we construct  a complete set of eigenfunctions of
the operator ${\bf H}$. They are parametrized by the magnetic quantum number $m$,
the momentum $p$ in the direction of the $x_{3}$-axis and the number $n$ of
an eigenvalue $\lambda_{m,n}(p)$ of the operator $H_m(p)$. Thus, if we
set
\bel{Psi1}
{\pmb \psi}_{n,m, p}(r,\theta, x_{3})=e^{ip x_{3}} e^{im\theta } \psi_{n, m}(r,p),
\ee
then
\bel{Psi2}
{\bf H}{\pmb \psi}_{n,m, p}=\lambda_{n, m}(p) {\pmb \psi}_{n,m, p}.
\ee
  
In Section 3, we show that for all $n \in {\mathbb N}$ and $m
\in {\mathbb Z}$ :
\begin{itemize}
    \item 
    Under very general assumptions $\lambda_{n,m}(p) \to \infty$ as $p \to \infty$ (Proposition \ref{t31}).
    \item If $b (r) \to 0$ as $r \to \infty$, then $\lambda_{n,m}(p) \to
0$ as $p \to -\infty$ (Proposition \ref{t32}).
    \item If $b (r)$ admits a finite positive limit $b_{0}$ as $r \to \infty$,
    then $\lambda_{n,m}(p) \to (2n-1)b_{0}$ for all $m$ as $p \to
-\infty$ (Proposition \ref{t33}).
    \item If
$b (r) \to \infty$ as $r \to \infty$, then $\lambda_{n,m}(p) \to
\infty$ as $p \to -\infty$ (Proposition \ref{t33}).
\end{itemize}

Related results concerning the dispersion curves  for
Schr\"odinger operator with constant magnetic fields defined on
unbounded domains $\Omega \subset \rd$ have been obtained in
\cite{GS} (the case where $\Omega$ is a strip) and in \cite{dBP},
\cite[Section 4.3]{H} (the case where $\Omega$ is a
half-plane). 

In Theorem \ref{f31} we formulate the main spectral results which
follow from the asymptotic properties of the dispersion curves
$\lambda_{n,m}(p)$, $p \in \re$. First, the analyticity and the
asymptotics as $p \to  \infty$ of $\lambda_{n,m}(p)$ imply
immediately that the spectra $ \sigma({\bf H}_m)$ and  $\sigma({\bf H})$  of the operators ${\bf H}_m$, $m \in
    {\mathbb Z}$,
and ${\bf H}$ are purely absolutely continuous.  Moreover, 
 \bel{kr37}
    \sigma({\bf H}_m)   = [{\mathcal
    E}_m,\infty),   \quad m \in {\mathbb
    Z},\q \sigma({\bf H})   = [{\mathcal
    E}_0,\infty),
    \ee
    where
    \bel{kr34}
    {\mathcal E_m}  =   \inf_{p \in \re} \lambda_{1,m}(p)\geq 0.
    \ee
Next, in the case
where the magnetic field tends to   $  0$ as
$r \to \infty$, the spectra of ${\bf H}_m$, $m \in {\mathbb Z}$,
  coincide with $[0,\infty)$ and have infinite
multiplicity. On the other hand, in the case where the magnetic
field tends as $r \to \infty$ to a positive finite limit, or to
infinity,   we have that
 ${\mathcal E}_m > 0$ for all
$m \in {\mathbb Z}$  and each of the spectra $\sigma({\bf H}_m)$  contains infinitely many thresholds. 

Further, in Section 4, we obtain a convenient formula
for the derivatives $\lambda_{n,m}'(p)$ which play the role of asymptotic group velocities.
 Our formula for $\lambda_{n,m}'(p)$ yields sufficient
conditions    (see Theorem \ref{t41}) for   positivity of these functions. The
leading example when these conditions are met, is
\bel{bb}
b(r)=b_{0}r^{-\delta} , \quad b_{0} >0, \quad \delta\in [0,1],
\ee
and  $m \neq 0$. If $\delta=1$, this result remains true for all $m$ (cf.   \cite{Y}). On the contrary, if $\delta =0$ and $m=0$, then $\lambda_{n,0}'(p) <0$ for all $n$
on  some  interval of  $p$ (lying on the negative half-axis). Similar results concerning for the dispersion
curves for the Schr\"odinger operator with constant magnetic
field, defined on the half-plane with Dirichlet (resp., Neumann)
boundary conditions, can be found in \cite{dBP} (resp., \cite{DH}
and \cite[Section 4.3]{H}).

  Finally, in Section 5 we discuss the long-time behaviour of a quantum particle.
    The time evolution of a quantum system  is determined by the unitary groups
$\exp{(-i{\bf H}_m t)}$, $m \in {\mathbb Z}$, so that an analysis of its asymptotics as $t\to\pm\infty$ relies on spectral properties of the operators ${\bf H}_m$.
Since these operators have discrete spectra, a quantum particle remains localized in the $(x_{1}, x_{2})$-plane. Its propagation in the $x_{3}$-direction is governed by the group velocities $\lambda_{n,m}'(p)$.  In particular, the condition $\lambda_{n,m}'(p)>0$ for all $n \in {\mathbb N}$ and $p \in {\mathbb R}$ implies that a   quantum particle with the magnetic quantum number $m$ propagates as $t\to+\infty$  in the positive direction of the $x_{3}$-axis.

Let us compare these results with  the long-time behaviour of a classical  particle in magnetic field \eqref{magnfi}. As shown in \cite{Y}, the function $x_{3}'(t)$ is periodic with  period $T$ determined by initial conditions. Its  drift $x_{3} (T) -  x_{3} (0)$ over the period is nonnegative if $b(r)\geq 0$ and $b'(r)\geq 0$. Morever, it is strictly positive if $b'(r)> 0$ for all $r$. In the   case $b(r) = {\rm const}$ it is still strictly positive if the angular momentum $m$ of a particle is not zero. Thus, our results
for functions \eqref{bb} correspond completely to the classical picture if $\delta=1 $ or   $\delta\in [0,1)$ and $m\neq 0$. In the case   $\delta =0$ and $m =0$ the behaviour of quantum and classical particles turn out to be qualitatively different.


\section{ Hamiltonians and their diagonalizations}
  \setcounter{equation}{0}
  
  
  Here we give precise definitions of the Hamiltonians and  discuss their reductions due to the cylindrical symmetry.
  
  \medskip
  
    {\bf 2.1.}
    For an arbitrary magnetic potential $A : \re^3 \to \re^3$ such that $A \in
  L^2_{\rm loc}(\re^3)^3$, the self-adjoint Schr\"odinger operator \eqref{MH} can be defined via its quadratic form
    \bel{w1}
{\bf h}[{\bf u}]=  \int_{\re^3} \left| i (\nabla {\bf u})(x)+ A(x){\bf u}(x)\right|^2dx, \q x=(x_{1}, x_{2}, x_{3}) .
  \ee
  It is easy to see that this form is closed on the set of functions  $u \in
  L^2 (\re^3) $ such that $\nabla {\bf u}   \in
  L^1_{\rm loc} (\re^3)^3$ and $i\nabla {\bf u} + A {\bf u} \in
  L^2 (\re^3)^3$. 
  Similarly, if $a \in L^2_{\rm loc}(\re^2)$, then the self-adjoint Schr\"odinger operator \eqref{u1ax} can be defined via its quadratic form
\bel{w1a}
h[u;p]=  \int_{\re^2} \left(\left| (\nabla u)(x)\right|^2 + (a(x)+p)^2 |u(x)|^2\right)dx, \q x=(x_{1}, x_{2}),   \quad  p \in \re.
  \ee
  This form is closed on the set of functions  $u \in
  L^2 (\re^2) $ such that integral \eqref{w1a} is finite. Clearly, this set does not depend on the parameter $p\in \re$. 
  
    Let ${\mathcal F}: L^2(\re^3    )\to L^2(\re; L^2(\re^2 )     )$ be the Fourier transform
with respect to $x_3$, i.e.
$$
({\mathcal F}{\bf u})(x_1,x_2,p) = \frac{1}{\sqrt{2\pi}} \int_{\re}
e^{-ix_3p} {\bf u} (x_1,x_2,x_3) dx_3. 
$$
  If $A(x)$ is given by formula \eqref{w3}, then 
  \[
  {\bf h}[{\bf u}]=\int_{-\infty}^\infty h[ ({\mathcal F} {\bf u})(p);p] dp
  \]
  which implies the equation
     \bel{w5}
  (  {\mathcal F} {\bf H} {\bf u})(x_{1}, x_{2}; p)=   (  H(p) {\mathcal F} {\bf u})(x_{1}, x_{2}; p).
    \ee
    This equation can be regarded as a ``working" definition of the operator  ${\bf H}$.
    
    \medskip
  
  {\bf 2.2.}
   Assume now   that the function $a$ in \eqref{w3} depends only
    on $r$, and
    \bel{kr8}
    a \in L^2_{\rm loc}([0,\infty);rdr).
    \ee
   If we separate variables in the cylindrical coordinates $(r,\theta, x_{3} )$ and denote by 
   ${\goth H}_m \subset L_2({\re}^3)$   the subspace of functions 
${\bf f}(r,x_{3})e^{im\theta }$ where ${\bf f} \in
L_2({\re}_+ \times {\re};rdr dx_{3})$   and  $m \in {\mathbb Z} $ is the magnetic quantum number, then  
\[
    L_2({\re}^3)=\bigoplus_{m \in {\mathbb Z}}  {\goth H}_m.
\]
The subspaces ${\goth H}_m$ are invariant with respect to $\bf H$ so that
  restrictions ${\bf H}_m$ of ${\bf H}$ on ${\goth H}_m$ are related with $\bf H$ by formula 
 \bel{Hm}
{\bf H}=\bigoplus_{m \in {\mathbb Z}}{\bf
H}_m.
\ee 
 Every ${\goth H}_m$
can obviously be identified with    the space $L^2(\re_+\times \re;rdr dx_{3})=:{\goth H}=$;
then   ${\bf H}_m$ is identified with   operator \eqref{Hm1}.

Quite similarly, if ${\cal H}_m
\subset L_2({\re}^2)$ is
the subspace of functions $f(r) e^{im \theta}$ where $f \in L_2({\re}_+;rdr)$, then
\[
    L_2({\re}^2)=\bigoplus_{m \in {\mathbb Z}} {\cal H}_m.
\]
The subspaces ${\cal H}_m$ are invariant with respect to $H(p)$ so that
  restrictions $H_m(p)$ of $H (p)$ on ${\cal H}_m$ are related with $H (p)$ by formula \begin{equation}\label{eq:HHm1}
H (p)=\bigoplus_{m \in {\mathbb Z}}   H_m(p).
\end{equation}
 Every ${\cal H}_m$
can obviously be identified with    the space $L^2(\re_+ ;rdr  )=: {\cal H}$;
then   $H_m(p)$ is identified with   operator \eqref{u1a}.

  Let ${\mathcal F}_{m} : {\goth H}_m \to L^2(\re; L^2(\re_+ ;rdr  )    )$ be the restriction of ${\mathcal F}$ on the subspace ${\goth H}_m$.
Then we have (cf. \eqref{w5})
     \bel{w5m}
   ( {\mathcal F}_{m} {\bf H}_{m} {\bf f}) (r; p) =
(  H_{m}(p){\mathcal F}_{m}  {\bf f})(r; p).
    \ee

Sometimes it is more convenient to consider instead of $H_{m}(p)$ the operator
\bel{Lm}
 L_m(p) = r^{1/2} H_m(p) r^{-1/2} =   -  \frac{d^2}{dr^2}    +
\frac{m^2-1/4}{r^2} + (a(r)+p)^2  
\ee
acting in the space $L^2(\re_+)$ and
       unitarily equivalent to the operator  $H_m(p)$. It is easy to see that the operator  $L_m(p)$ corresponds to the quadratic form
         \bel{w1ay}
    l_{m}[g; p] = \int_0^{\infty} \left( |g'(r)|^2 + (m^2-1/4) r^{-2} |g (r)|^2 +
    (a(r)+p)^2 |g (r) |^2\right) dr,
    \ee
    defined originally on $C_0^{\infty}(\re_+)$, and then closed in
    $L^2(\re_+)$.
    
    \medskip

{\bf 2.3.} 
If
     \bel{j1}
     a(r) \to \infty \quad  {\rm as} \quad r \to \infty,
    \ee
then the spectrum of the operator $H_m(p)$, $p \in \re$, $m \in
{\mathbb Z}$, is discrete.  Thus, it consists of the increasing sequence
$\lambda_{n,m}(p)$ of simple eigenvalues.
Since $H_m(p)$,
$p \in \re $, is a Kato analytic family of type (B) (see
\cite[Chapter VII, Section 4]{K}), all the eigenvalues
$\lambda_{n,m}(p)$ are real analytic
functions of $p \in \re$. 
Moreover, $\lambda_{n,m}(p)>0$ because form \eqref{w1a} is strictly positive.

In view of formula \eqref{w5m}
spectral analysis of the operators $ {\bf H}_m$ reduces to a study of a family of functions
$\lambda_{n,m}(p) $, $n\in {\mathbb N}$. Indeed,
let $\Lambda_{n,m} $ be the operator   of multiplication by the function $\lambda_{n,m}(p) $ in the space $L^2 (\re)$. 
We denote by $\psi_{n,m}(r;p)$ real normalized  eigenfunctions (defined up to signs) of the operators $H_{m}(p)$ and introduce an isometric mapping 
\[  
   \Psi_{n, m} : L^2(\re    ) \rightarrow L^2({\re}_+\times {\re}; rdr d p) 
    \]
by the formula   
  \bel{Psi}
(\Psi_{n,m} w) (p)=  \psi_{n,m}(r,p) w(p).
\ee
 Then 
 \[
 L^2({\re}_+\times {\re}; rdr d p)=\bigoplus_{n\in {\Bbb N}} {\ran }\Psi_{n,m}
 \]
 and 
   \bel{UI}
      {\bf H}_{m}     =\bigoplus_{n\in {\Bbb N}} 
    {\cal F}_{m}^*  \Psi _{n,m} \Lambda_{n,m} \Psi_{n,m}^* {\cal F}_{m}. 
\ee
    Together with \eqref{Hm},
formulas \eqref{Psi} and  \eqref{UI} justify equations \eqref{Psi2} for functions \eqref{Psi1}.


    \section{Dispersion curves and spectral analysis }
    \setcounter{equation}{0}


{\bf 3.1.} 
In this subsection we consider the operators $H(p)$ acting in the space 
$    L^{2} (\rd)$ by formula \eqref{u1ax}. Under the assumption $  a \in L^{2}_{\rm loc}(\rd)$ they are correctly defined by their quadratic forms  \eqref{w1a}. If 
\bel{Di}
a(x )\to\infty\q {\rm as} \q |x|\to\infty,\q x=(x_{1}, x_{2}) , 
\ee
then the spectrum of $H(p)$ consists of eigenvalues $\lambda_{n }(p)$, $n \in {\mathbb N}$. We enumerate them in the increasing order with multiplicity taken into account. Our goal is to
investigate the asymptotic
behaviour of the eigenvalues $\lambda_{n}(p)$ as $p \to \infty$.
Below we denote by $C$ and $c$ different positive constants whose precise values are of no importance.

We use the following elementary

\begin{lemma} \label{l31}
Let $v(x)\geq 0$.
For an arbitrary $\varepsilon>0$, we have the inequality
\bel{2y}
\int_{\rd} v(x) |  u(x)|^2 dx \leq C \sup_{x\in {\rd}}\left( \int_{|x-y|\leq \varepsilon} v^2(y)dy\right)^{1/2}
 \int_{\rd}  \left( \varepsilon |\nabla u(x)|^2   +\varepsilon^{-1}   |  u(x)|^2\right) dx 
\ee
provided the supremum in the right-hand side is finite.
    \end{lemma}

\begin{proof}
Let $\Pi_{\varepsilon}\subset \rd$ be a square of length $\varepsilon$. We  proceed from the estimate
  \[
 \left(   \int_{\Pi_{\varepsilon}}   |  u(x)|^4 dx \right)^{1/2}\leq C
  \left(   \varepsilon \int_{\Pi_{\varepsilon}}    |\nabla u(x)|^2 dx 
  + \varepsilon^{-1} \int_{\Pi_{\varepsilon}}    |  u(x)|^2 dx\right)
    \] 
which follows from the Sobolev embedding theorem by a scaling transformation.
Using the Schwarz inequality, we deduce from this estimate that
 \bel{2x}
    \int_{\Pi_{\varepsilon}} v(x)  |  u(x)|^2 dx  \leq C
  \left(    \int_{\Pi_{\varepsilon}} v^2(x)   dx \right)^{1/2} 
  \left(   \varepsilon \int_{\Pi_{\varepsilon}}    |\nabla u(x)|^2 dx 
  + \varepsilon^{-1} \int_{\Pi_{\varepsilon}}    |  u(x)|^2 dx\right).
    \ee 
Let us split the space $\rd$ in  the lattice of squares $\Pi_{\varepsilon}^{(n)}$ of length $\varepsilon$. Applying \eqref{2x} to every $\Pi_{\varepsilon}^{(n)}$ and summing over all
$n$, we  arrive at \eqref{2y}.
\end{proof}

In the following assertion we do not assume \eqref{Di}.
   
\begin{pr} \label{l31xx}
Let  $  a \in L^{2}_{\rm loc}(\rd)$.
Set $a_{-}(x)=\max\{-a(x),0\}$,
 \[
\alpha(\varepsilon) = \sup_{x\in {\rd}}  \int_{|x-y|\leq \varepsilon} a_-^2(y)dy 
     \]
 and suppose that $\alpha(\varepsilon)\to 0$  as $\varepsilon\to 0$.  Then we have
    \bel{2}
    \liminf_{p \to \infty} p^{-2} \inf \sigma(H(p))
    \geq 1.
    \ee
\end{pr}

\begin{proof}
Applying estimate  \eqref{2y}
with $\varepsilon=p^{-1}$  to the function $v=a_{-}$, we find that
 \begin{eqnarray*}
\int_{\rd} \left(|\nabla u|^2 +  (p + a)^2  \right) |u|^2   dx \geq \int_{\rd}
\left( |\nabla u|^2 + (- 2p a_-  +p^2) |u|^2 \right) dx
\\
\geq
\int_{\rd} \left(|\nabla u|^2 +  p  ^2  \right) |u|^2   dx
-C \sqrt{\alpha(p^{-1}) }  
\int_{\rd} \left(|\nabla u|^2 +  p^2  \right) |u|^2   dx.
     \end{eqnarray*}
 Since $\alpha(p^{-1}) \to 0$ as $p\to\infty$, this implies  \eqref{2}.
  \end{proof}

\begin{pr} \label{t31} 
Let $  a \in L^{2}_{\rm loc}(\rd)$ and let condition \eqref{Di} be satisfied. 
 Then, for all  $n \in {\mathbb N}$,  we have
    \bel{kr4}
    \lambda_{n}(p) = p^2(1 + o(1)), \quad p \to \infty.
    \ee
\end{pr}

\begin{proof}
Under condition \eqref{Di} the function $a_{-}$ has compact support  so that
 we can use Proposition~\ref{l31xx} and  estimate \eqref{2} implies
    \bel{kr4a}
    \liminf_{p \to \infty} p^{-2} \lambda_{n }(p) \geq 1.
    \ee
    Set $G (\varepsilon)=-\Delta + (1 + \varepsilon^{-1})
a^2(x) $, $\varepsilon > 0$. The spectrum of $G (\varepsilon)$ is discrete;
let $ \nu_{n }$, $n \in {\mathbb N}$,  be
the increasing sequence of its eigenvalues. By the elementary
inequality 
$$(a+p)^2 \leq (1 + \varepsilon^{-1}) a^2 + (1 +
\varepsilon) p^2, \q \varepsilon > 0, 
$$
we have $H (p) \leq
G (\varepsilon) + (1 + \varepsilon) p^2$ so that by the minimax
principle   
$$\lambda_{n }(p) \leq \nu_{n }(\varepsilon) +
(1 + \varepsilon) p^2.
$$
Therefore,  for all $\varepsilon>0$,
$$\limsup_{p \to \infty} p^{-2}
\lambda_{n }(p) \leq 1 + \varepsilon,  
$$
which combined with \eqref{kr4a} yields \eqref{kr4}.
\end{proof}

\begin{follow} \label{t31f} 
Suppose that the function $a$ depends on $r$ only.
Let conditions \eqref{kr8} and \eqref{j1} be satisfied. 
 Then, for all  $n \in {\mathbb N}$, $m \in {\mathbb
Z}$,  we have
  \[
    \lambda_{n,m}(p) = p^2(1 + o(1)), \quad p \to \infty.
    \]
\end{follow}

\medskip

    {\bf 3.2.}
    {\it From now on we always assume that the function $a$ depends on $r$ only  and that conditions \eqref{kr8} and \eqref{j1}  are satisfied}.  
   In this subsection we investigate the asymptotics
as $p \to -\infty$ of the eigenvalues $\lambda_{n,m}(p)$ of the operators $H_{m}(p)$.
Actually, it is more convenient to work with the operators $L_{m}(p)$ acting in the space 
$    L^{2} (\re_{+})$ by formula \eqref{Lm}.
We suppose that the function $a$  is differentiable at least for suffiently big $r$ and  formulate the results in terms of the  function $b(r)=a'(r)$ related to the magnetic field by formula \eqref{magnfi}.

Remark first that if     $k= -p>0$ is big enough,  then the equation
    \bel{e5}
    a(r) = k
    \ee
  has at least one solution.
We denote by $\rho_k$ the greatest solution of (\ref{e5}). 
Clearly, $\rho_k\to\infty$ as $k\to\infty$.

\begin{pr} \label{t32}
   Suppose  that  
    \bel{e3}
    \lim_{r\to\infty} b(r) = 0.
    \ee
Then for each $n \in {\mathbb N}$ and $m \in {\mathbb Z}$ we have
    \bel{e4}
    \lim_{k \to\infty} \lambda_{n,m}(-k) = 0.
    \ee
\end{pr}

\begin{proof}
Set
 \bel{E} 
{\bf b}(r) = \sup_{x \geq r} |b (x)| \q {\rm and}\q  \gamma_k=
{\bf b} (\rho_k)^{-1/2}.
\ee
 Let us fix $n \in
{\mathbb N}$.  We pick a   function $\phi_1 \in
C_0^{\infty}(\re)$  such that ${\rm supp}\,\phi_1 = \left[0, \frac{1}{2n}\right]$ and, for $n>1$, set
$$
\phi_j(x)  = \phi_1(x - (j-1)/n), \quad x \in \re, \quad j=2,\ldots,
n.
$$
For $k>0$ large enough, we put
 \bel{e5xx}
\varphi_{j}(r;k) = \gamma_k^{-1/2}
\phi_j \left(\frac{r-\rho_k}{\gamma_k}\right), \quad r \geq 0,
\quad j=1,\ldots, n.
\ee

We will prove now that for quadratic form \eqref{w1ay}
    \bel{e8}
    \lim_{k \to \infty} l_{m} [\varphi_{j}(k); -k] = 0.
    \ee
 It follows from \eqref{e5xx} that
\bel{e9}
    \int_0^{\infty} |\varphi'_{j} (r;k)|^2 dr \leq C  \gamma_k^{-2}
    \ee
with $C$ independent of $k$. Further, since $\supp  \varphi_{j}(k)  
\subset [\rho_k, \rho_k + \gamma_k ]$, we have
    \bel{e10}
    \int_0^{\infty} r^{-2}|\varphi_{j} (r;k)|^2 dr \leq C    \rho_k ^{-2}.
    \ee
  Similarly,
    \bel{e11}
    \int_0^{\infty} (a(r)-k)^2 |\varphi_{j} (r;k)|^2 dr
    \leq C    \sup_{r \in (\rho_k, \rho_k +  \gamma_k )} (a(r)-k)^2.
    \ee
    Using the condition $a(\rho_{k})=k$, we obtain, for $r\geq \rho_{k}$, the bound
    \[
  (a(r)-k)^2 =   (a(r)- a(\rho_{k}))^2 = \left( \int_{\rho_{k}}^r b (s)ds\right)^2\leq (r-\rho_{k})^2 {\bf b}^2(\rho_k)   
    \]
    where ${\bf b} $ is function (\ref{E}). Thus, the right-hand side in (\ref{e11}) is bounded by $C \gamma_k^2 {\bf b}^2(\rho_k)$. Putting together this result with
     inequalities (\ref{e9}), (\ref{e10}) and taking into account (\ref{E}),  we
get
$$
l_{m}[\varphi_{j}(k); -k] \leq C \left({\bf b} (\rho_k) +  
\rho_k^{-2}  \right).
$$
This yields (\ref{e8}).

Let us use now that the supports of the functions $\varphi_{j}(k)$,
$j=1,\ldots, n $, are disjoint and set
 \bel{j9a}
    {\mathcal L}_{n}(k) = {\rm span}\, \left\{\varphi_{1}(k), \ldots, \varphi_{n}(k)\right\}.
    \ee
    Then $ {\rm dim}\;  {\mathcal L}_{n}(k) = n$ and according to (\ref{e8})
    $  l_{ m} [\varphi (k); -k] \to 0$ as $k \to \infty$ 
for all $\varphi (k)\in {\mathcal L}_{n}(k) $ with $\|\varphi (k)\|=1$. By the mini-max principle this  implies (\ref{e4}). 
\end{proof}

The proof of Proposition~\ref{t33} relies on a comparison of the operator $L_{m}(-k)$ with  the ``model" operator
  \bel{j2c}
        T(k) = -\frac{d^2}{dx^2} +
    b^2 (\varrho_k)  (x-\varrho_k)^2, \quad x \in \re,
    \ee
    acting in the space $L^2(\re)$. Let $f_j$ be
the normalized in $L^2(\re)$ real-valued eigenfunctions (defined up to sign) of the
harmonic oscillator, i.e.
   \bel{HO}
-f_j''(x) + x^2 f_j(x) = (2j-1) f_j(x), \quad x \in \re, \quad j
\in {\mathbb N}.
 \ee
 Then
  \bel{j3}
    \psi_{j}(x;k)  =
    b (\varrho_k)^{1/4} f_j (b (\varrho_k)^{1/2}(x-\varrho_k)) 
        \ee
    are normalized eigenfunctions of the operator $T (k)$, that is
       \bel{w30}
    T (k) \psi_{j} (k) = b (\varrho_k) (2j-1) \psi_{j} (k), \quad j
    \in {\mathbb N}.
    \ee
   
 The proof of the following result   follows the general lines of the
proof of \cite[Theorem 11.1]{CFKS}.

\begin{pr} \label{t33}
Suppose that $a(r)$ is locally semibounded from above. For $r>0$ large enough, we assume that the function $b(r)$ is   differentiable and that conditions 
 \bel{j2}
    b (r) > 0, 
    \ee
     \bel{j2bX}  
  \lim_{r \to\infty}  r^{2} b (r) =\infty,
    \ee 
     as well as 
     \bel{j3a}
    \lim_{r \to \infty}  b (r)^{-3}  {\bf b}_{1}^2  (r) = 0, \q {\rm where} \q {\bf b}_{1}(r)   = \sup_{r/2\leq x \leq 3r/2} |b'(x)|,
    \ee
are satisfied. Let also 
  \bel{j2b}  
  \lim_{k \to\infty}  k^{-2} b (\rho_{k}) =0. 
    \ee  
 Then, for all  $n \in {\mathbb N}$, $m \in {\mathbb Z}$, we have
    \bel{j4}
    \lambda_{n,m}(-k) = b (\varrho_k) (2n-1 + o(1)), \quad k \to
    \infty.
    \ee
\end{pr}

\begin{proof}
  Due to the minimax principle,
it suffices to show that:\\
(i) For each $n \in {\mathbb N}$ and sufficiently large $k$ there exists a subspace
${\mathcal L}_{n} (k)$ of $L^2(\re_+)$ such that ${\rm
dim}\,{\mathcal L}_{n} (k) = n$, ${\mathcal L}_{n} (k) \subset
D(L_m(-k))$, and for each $\varphi (k) \in {\mathcal L}_{n} (k)$   we have
    \bel{j5}
    \langle L_{m}(-k) \varphi (k), \varphi (k)\rangle  \leq b (\varrho_k) (2n-1 + o(1))\|\varphi (k) \|^2,
    \quad k \to \infty.
    \ee
(ii) For each $n \in {\mathbb N}$ there exists a bounded operator
$R_n (k)$ such that ${\rm rank}\,R_n (k)\leq n-1$ (hence, $R_1 (k) = 0$), and
    \bel{j6} L_m(-k)  \geq b (\varrho_k) (2n-1 + o(1))I + R_n (k),
    \quad k \to \infty.
    \ee
    
 We pick
$\gamma_k > 0$ such that
    \bel{j7}
    \gamma_k \to 0,
    \ee
    \bel{j8}
    \gamma_k \, \varrho_k b(\varrho_k)^{1/2} \to \infty,
    \ee
    \bel{j8a}
    \gamma_k^{-3} \, b(\varrho_k)^{-3/2}{\bf b}_{1} (\varrho_k)\to 0
    \ee
    as $k\to\infty$.
Note that \eqref{j8} is compatible with \eqref{j7} due to
\eqref{j2bX}, and \eqref{j8a} is compatible with \eqref{j7} due to
\eqref{j3a}.

\medskip

{\em Proof} of (i). 
 Let $\zeta \in C_0^{\infty}(\re)$ be such that $0\leq \zeta(x) \leq 1$ , $\zeta(x) = 1$ for $|x|\leq 1/2$ 
and $\supp \zeta = [-1,1]$. 
For $k$ large enough, set
 \bel{j9b}
    \zeta (r;k) = \zeta(\gamma_k  b (\varrho_k)^{1/2}(r-\varrho_k)), \quad r \in \re_{+}, 
    \ee
    and
    \bel{j9}
     \varphi_{j}(r;k)  = \psi_{j}(r;k) \zeta (r;k),
    \quad r \in \re_{+},
    \quad j \in {\mathbb N},
    \ee
the functions $\psi_{j} (r;k)$ being defined in \eqref{j3}.
 It follows from \eqref{j8} that
$$
\supp {\varphi_{j}} (k) =  [ \varrho_k -  \gamma_k^{-1}
b (\varrho_k)^{-1/2}, \varrho_k +  \gamma_k^{-1}b (\varrho_k)^{-1/2} ] \subset [\varrho_k/2, 3\varrho_k/2] 
$$
  and, in particular,  ${\varphi_{j}}(k) \in D(L_m(-k))$.
Note that
    \bel{j19a}
    \langle \varphi_{j}(k), \varphi_{l}(k)\rangle_{L^2(\re_+)} =
  \delta_{jl} -
    \int_{\re} \psi_{j}(x;k)\psi_{l}(x ; k) (1 - \zeta^2(x; k))dx =
    \delta_{jl} + o(1)   
    \ee
    as $k \to \infty$.  Indeed, the integral here can be estimated by
$$
  \int_{\re} |f_{j}(x)f_{l}(x)| (1 -
\zeta^2(\gamma_k x) )dx \leq
\int_{|x| \geq (2\gamma_k)^{-1} } |f_{j}(x)f_{l}(x)| dx 
$$
which tends to zero according   to \eqref{j7}. In particular, \eqref{j19a} implies that for all
$n \in {\mathbb N}$  the functions $\varphi_{1}(k), \ldots, \varphi_{n}(k)$ are
linearly independent if $k$ is large enough. Thus, the space
 $   {\mathcal L}_{n}(k)$ defined by  \eqref{j9a} has dimension $n$.
 
Let us set
 \bel{psi}
 \psi (x;k )   = \sum_{j=1}^n c_j \psi_j(x;k ), \q c_j \in {\mathbb C},
 \ee
  $  \varphi (r;k)  = \psi (r;k) \zeta (r;k)$ and consider $ \langle L_{m}(-k) \varphi (k), \varphi (k)\rangle$. 
 Integrating by parts, we find that
\[
-2{\rm Re}  \, \langle  \psi'(k) \zeta'(k) , \psi (k)\zeta(k)\rangle-\langle  \psi(k) \zeta''(k) , \psi (k)\zeta(k)\rangle = \| \psi(k) \zeta'(k)\|^2
\]
so that  
 \begin{eqnarray} 
    \langle L_{m}(-k) \varphi (k), \varphi (k)\rangle  &=&
{\rm Re} \,  \langle  -\psi''(k) + (a(r)-k)^2 \psi (k), \psi (k)\zeta^2(k)\rangle 
\nonumber\\
&+&
        \| \psi(k) \zeta'(k)\|^2
    + (m^2-1/4) \| r^{-1} \varphi (k)  \|^2.
\label{j5x1}   \end{eqnarray}
We assume that $\| \varphi (k)\|=1$ and hence according to \eqref{j19a}
$\| \psi(k)\|=1+ o(1)$.
The second and third terms in the right-hand side of \eqref{j5x1} are negligible.
 Indeed, differentiating \eqref{j9b} and using condition \eqref{j7}, we find that
 \bel{j12}
    \|\psi  (k)\zeta' (k)\|^2 = O(b (\varrho_k) \gamma_k^2) = o(b (\varrho_k)).
    \ee
 Since $ r^{-1} \leq 2 \varrho_k^{-1}$ on the support of $\varphi  (k)$, relation \eqref{j2bX} implies
    \bel{j11}
    \|r^{-1} \varphi(k)   \|^2= O (\varrho_k^{-2})  = o( b (\varrho_k)),
    \quad k \to \infty.
    \ee
 
 Further we consider the first term in the right-hand side of \eqref{j5x1}.
 It follows from equation \eqref{w30} that
    \bel{j14}
    -\psi_{j}''(k) + (a(r)-k)^2 \psi _{j}(k) =
    b (\varrho_k)(2j-1) \psi_{j} (k) +
 \alpha(k) \psi_{j}(k)
    \ee
    where the function
    \bel{alpha}
 \alpha(r;k)=  (a(r)-k)^2 - b ^2 (\varrho_k)     (r-\varrho_k)^2. 
    \ee
    Let us estimate the right-hand side. In view of the equation $a(\rho_{k})=k$,  a second-order Taylor expansion of $a$ at
$\varrho_k$ yields
$$
a(r) = k + b (\varrho_k) (r-\varrho_k) + \int_{\varrho_k}^r b '(s)
(r-s) ds.
$$
Therefore,
$$
 \alpha(r;k) = 2 b(\varrho_k) (r-\varrho_k)\int_{\varrho_k}^r
b'(s) (r-s) ds + \left(\int_{\varrho_k}^r b '(s) (r-s)
ds\right)^2,
$$
and hence
$$
  |\alpha(r;k)| \leq b (\varrho_k) {\bf b}_{1}(\varrho_k)
| r-\varrho_k |^3 + 4^{-1}{\bf b}_{1}^2 (\varrho_k) (r-\varrho_k)^4 
$$
\[
   \leq  \gamma_k^{-3} b (\varrho_k)^{-1/2} {\bf b}_{1}(\varrho_k) +
   4^{-1} \gamma_k^{-4} b (\varrho_k)^{-2} {\bf b}_{1}^2(\varrho_k),
    \]
 provided that $|r-\varrho_k| \leq \gamma_k^{-1} b (\varrho_k)^{-1/2}$.
 In view of conditions \eqref{j7} and \eqref{j8a}, this gives us the estimate
    \bel{j17n}
   \sup_{|r-\varrho_k| \leq \gamma_k^{-1} b (\varrho_k)^{-1/2}} |\alpha(r;k) | =o(b (\rho_{k}))
    \ee
 so that
  \[
  \left( -\psi_{j}''(k) + (a(r)-k)^2 \psi _{j}(k) \right) \zeta(k)=
    b (\varrho_k)(2j-1) \varphi_{j} (k)+  o(b (\rho_{k})).
    \]
  Thus,  using also \eqref{j19a} we obtain that 
     \begin{eqnarray} 
 {\rm Re} \,  \langle  -\psi''(k) + (a(r)-k)^2 \psi (k), \psi (k)\zeta^2(k)\rangle   &=&
b (\varrho_k) \sum_{j,l=1}^n (2j-1) c_{j}\bar{c}_{l}  \langle  \varphi _{j}(k), \varphi_{l} (k)\rangle +  o(b (\rho_{k}))
\nonumber\\
&\leq&
    b (\varrho_k) (2n-1)   +  o(b (\rho_{k})).
\label{j5x2}   \end{eqnarray}
     Together with \eqref{j12} and \eqref{j11}, this implies estimate
\eqref{j5}   for each $\varphi (k) \in {\mathcal L}_{n}(k)$.

\medskip

{\em Proof} of (ii).
 Let functions $\zeta \in
C_{0}^{\infty}(\re)$ and  $\eta
\in C^{\infty}(\re)$ satisfy
$
\zeta^2(x)  + \eta^2(x)  = 1$, $ x \in \re$;
moreover, as before,  we require that $0\leq \zeta(x) \leq 1$ , $\zeta(x) = 1$ for $|x|\leq 1/2$ 
and $ \supp \zeta = [-1,1]$. 
By analogy with \eqref{j9b} set 
\bel{j21q}
\eta(r;k)  = \eta(\gamma_k
b (\varrho_k)^{1/2} (r-\varrho_k)), \q r \in \re_{+}. 
\ee
Then we have
  \[
    \zeta^2(r;k) + \eta^2(r;k) = 1, \quad r \in \re_{+}.
    \]
We proceed from the    localization formula (known as the IMS   formula -- see e.g. \cite[Section 3.1]{CFKS})  
  \[
    L_m(-k) = \zeta(k) L_m(-k)\zeta (k)  + \eta (k) 
    L_m(-k)\eta (k)  - \zeta'(k)^2 -  \eta'(k)^2,
    \]
    where $\zeta(k)$, $\eta (k)$, $\zeta'(k)$ and $\eta'(k)$ are understood as operators of multiplication by the functions $\zeta(r,k)$, $\eta (r,k)$, $\zeta'(r,k)$ and $\eta'(r,k)$, respectively.
   According to \eqref{j7} it follows from definitions  \eqref{j9b} and \eqref{j21q} that
  \bel{j22}   
\max_{r\in\re_{+}}\,(\zeta' (r,k) ^2 +  \eta'(r,k) ^2 )= O\left(\gamma_{k}^2 b (\varrho_k)\right)
= o\left(b (\varrho_k)\right), \quad k
\to \infty. 
\ee 

Next, we check that
 \bel{j23w}
       \eta (k) L_{m}(-k)  \eta (k)\geq \nu_k b (\varrho_k) \eta^2 (k)
    \ee
with $\nu_k \to \infty$ as $k \to \infty$. By virtue of the
 Hardy inequality 
  \[
    \eta (k) \left( -\frac{d^2}{dr^2} +
    \frac{m^2-1/4}{r^2}\right) \eta (k) \geq 0,
    \]
    it suffices to check that
    \bel{j24}
    (a(r) - k)^2 \geq \nu_k b(\varrho_k) 
    \ee
    for 
 \bel{j24c}
    r  \geq  \varrho_k+ \left(2\gamma_k
b (\varrho_k)^{1/2}\right)^{-1}=:\varrho_k^{(+)}\q {\rm and} \q     r  \leq  \varrho_k-\left(2\gamma_k
b (\varrho_k)^{1/2}\right)^{-1}=:\varrho_k^{(-)}.
    \ee
   According to \eqref{j2} there exists $r_{0}$ such that the function $a(r)$ is increasing for $r\geq r_{0}$.
   
    Let first $r\geq r_{0}$.
      Then
    \bel{j24d}
   | a(r) - k | = | a(r) -  a(\varrho_k)| \geq \pm (a(\varrho_k^{(\pm)}) -  a(\varrho_k))
   \ee
   if $\pm (r- \varrho_k^{(\pm)})\geq 0$ and  $r\geq r_{0}$.  It follows from definition 
   \eqref{j24c} of the numbers $\varrho_k^{(\pm)} $ that
     \bel{j24e}
    a(\varrho_k^{(\pm)}) -  a(\varrho_k) = \int_{\varrho_k}^{\varrho_k^{(\pm)}}b (s)ds= 
    \pm (2\gamma_k)^{-1} b (\varrho_k)^{1/2}+\int_{\varrho_k}^{\varrho_k^{(\pm)}}(b (s)-b (\varrho_k))ds.
   \ee
The absolute value of the integral in the right-hand side can be estimated by   
\[
 \left|    \int_{\varrho_k}^{\varrho_k^{(\pm)}}ds  \int_{\varrho_k}^s | b '(\sigma)|d\sigma\right|
     \leq  2^{-1} {\bf b}_{1} (\rho_{k}) (\varrho_k^{(\pm)}-\varrho_k)^2 =
     8^{-1} {\bf b}_{1} (\rho_{k})   \gamma_k^{-2}b (\varrho_k)^{-1}
   \]  
   where the function ${\bf b}_{1} $ is defined in \eqref{j3a}. By virtue of conditions \eqref{j7} and \eqref{j8a} this expression is $o(\gamma_k^{-1} b (\varrho_k)^{1/2})$ as $k\to\infty$. Therefore the absolute value of expression \eqref{j24e} is bounded from below by 
   $ (3\gamma_{k})^{-1} b(\varrho_k)^{1/2}$. Thus,   for $r\geq r_{0}$, estimate \eqref{j24}
   with $\nu_{k}=(3\gamma_{k})^{-2}\to\infty$   is a consequence of \eqref{j24d}. 
   
   If  $r\leq r_{0}$, we take into account that $a(r)$ is semibounded from above so that $(a(r) - k)^2\geq 2^{-1} k^2$. Hence estimate \eqref{j24d} with $\nu_{k}= 2^{-1} k ^{2}b (\varrho_k)^{-1} \to\infty$  is satisfied according to condition \eqref{j2b}.
   
 Putting together  definitions \eqref{Lm} and \eqref{j2c} of the operators $L_m(-k)$ and $T(k)$, we see that
 \bel{j27}
    \zeta (k) L_m(-k) \zeta (k)=
    \zeta (k) T (k) \zeta (k) +    \alpha(k) \zeta^2 (k),
    \ee
    where $\alpha(k)$ is the operator of multiplication by function \eqref{alpha}. The first term in the right-hand side
    is bounded from below by $  b (\varrho_k)   \zeta^2 (k)$
because  $b (\varrho_k)$ is the first eigenvalue of the operator
$T(k)$.  By virtue of \eqref{j17n} 
the second term satisfies the estimate
 \bel{j27Y}
\|  \alpha(k) \zeta^2 (k)\|=o(b (\varrho_k)) .
    \ee
      It follows that operator \eqref{j27}
  is bounded from below by $  b (\varrho_k)   \zeta^2 (k) -o(b (\varrho_k))I$.   Combining this result with  \eqref{j22} and \eqref{j23w}, we get
estimate  \eqref{j6} in the case $n=1$. 

If $n \geq 2$, we denote by $ P_{n}(k)$ the orthogonal
projection onto the span of the first $n-1$ eigenfunctions of the
operator $T (k)$. Then  $T(k)(I-P_{n}(k))\geq (2n-1)(I-P_{n}(k))$ and hence
$$
\zeta (k) T (k) \zeta (k) = \zeta (k) T (k)(I-P_n(k)) \zeta (k)
+ \zeta(k)T (k) P_n (k) \zeta (k)  
$$
    \bel{j28}
\geq  b (\varrho_k) (2n-1)
\zeta (k) (I-P_n (k)) \zeta (k) + \zeta (k) T (k) P_n (k) \zeta (k)   =
 b (\varrho_k) (2n-1) \zeta ^2 (k) + R_n (k)
    \ee
where
$$
R_n (k) = \zeta (k)( T (k) -  b (\varrho_k) (2n-1) I) P_n(k)
\zeta(k).
$$
Clearly, ${\rm rank}\,R_n (k)\leq n-1$.  Putting together \eqref{j22}, \eqref{j23w} and  \eqref{j27} --  \eqref{j28}, we obtain \eqref{j6}
in the case $n \geq 2$.  
\end{proof}

\begin{example} \label{ex}
Let $b(r)=b_{0}r^{-\delta}$,    $b_{0}>0$, $\delta \leq 1$, for sufficiently large $r$. Then 
${\bf b}_{1}(r)=b_{0}\delta r^{-\delta-1} $ and   conditions \eqref{j2} -- \eqref{j3a} are satisfied. Moreover, $\rho_{k}= c_{1} k^\nu$ and  $k^2 b (\rho_{k}) = c_2 k^{ -1-\nu}$ where $\nu= (1-\delta)^{-1}$ and $c_{1}, c_{2}>0$ if $\delta < 1$. If $\delta = 1$, then 
$\rho_{k}= \exp (b_{0}^{-1} k)$ and $k^2 b (\rho_{k}) = k^2 b_0 \exp (- b_{0}^{-1} k)$. In both cases condition \eqref{j2b} is also satisfied. Thus, Proposition~\ref{t33} implies the following results. 
If $\delta> 0$, then  $\lambda_{n,m}(p)\to 0$  as $p\to -\infty$ (this result follows also from Proposition~\ref{t32}).
If $\delta= 0$, then the functions $\lambda_{n,m}(p)$ have finite limits $b_{0}(2n-1)$ as $p\to -\infty$. If $\delta < 0$, then these functions tend to $+\infty$ as $p\to -\infty$.
\end{example}


\medskip

    {\bf 3.3.}
 Let us return to the Hamiltonians ${\bf H}_{m} $ and ${\bf H} $ defined in Section~2.

\begin{theorem} \label{f31}
Assume \eqref{kr8} and \eqref{j1}.

$(i)$ Then all operators ${\bf H}_m$, $m \in {\mathbb Z}$, and hence ${\bf H} $ are absolutely continuous and their spectra coincide with the half-axes defined by equations
 \eqref{kr37} and  \eqref{kr34}.

$(ii)$ If the hypotheses of Proposition~$\ref{t32}$ hold true, then ${\mathcal
    E}_m=0$ for all $m \in {\mathbb Z}$.
Moreover, the multiplicities of all spectra $\sigma({\bf H}_m)$ and hence of $\sigma({\bf H})$  are infinite.

   $(iii)$  Let the hypotheses of Proposition~$\ref{t33}$ hold true. If $b (r) \to \infty$, then
    the infimum in \eqref{kr34} is attained $($at a finite point$)$ so that for all $m \in {\mathbb Z}$
\[
   {\mathcal E_m} =\min_{p\in \re}\lambda_{1,m}(p) > 0.
    \]
    
$(iv)$ Let the hypotheses of Proposition~$\ref{t33}$ hold true. If
$b (r)$ admits a finite positive limit $b_{0}$ as $r \to \infty$, then 
${\mathcal E}_m \in (0,b_{0}]$ for all $m \in {\mathbb Z}$.
\end{theorem}

\begin{proof}
It suffices to prove only the assertions concerning the
operators ${\bf H}_m$. In view of decomposition  \eqref{UI} they reduce to corresponding statements about the operators $\Lambda_{n,m}$. These operators
are absolutely continuous because the eigenvalues
$\lambda_{n,m}(p)$ are real analytic   functions of $p
\in \re$ which are non constants  since according to Corollary~\ref{t31f}
$\lambda_{n,m}(p)\to\infty$ as $p\to\infty$. Moreover, we have that
\bel{kr34x}
  \sigma(\Lambda_{n,m})=[  {\mathcal E}_{n,m},\infty)\q {\rm where} \q   {\mathcal E}_{n,m} =  \inf_{p \in \re} \lambda_{n,m}(p)\geq 0
    \ee
   because $\lambda_{n,m}(p)>0$ for all  $p\in \re$.
    This implies relations   \eqref{kr37}   with ${\mathcal E_m}  $ defined by
     \eqref{kr34}.

In case (ii) it suffices to use that according to \eqref{e4}   ${\mathcal E}_{n,m}  =0$ and hence   $\sigma(\Lambda_{n,m})=[0,\infty)$ for all $m$ and $n$.

In case (iii)  Proposition~\ref{t33} implies that  $ \lambda_{n,m}(p)\to\infty$ as $p\to-\infty$ for all $n$ and $m$ so that 
\bel{kr34xy}
    {\mathcal E}_{n,m} =  \min_{p \in \re} \lambda_{n,m}(p) > 0
    \ee
and hence infimum   in \eqref{kr34} 
can be replaced by  minimum. 

In case (iv) we use that according to \eqref{j4} ${\mathcal E}_{n,m} \leq (2n-1)b_{0}$. Moreover, ${\mathcal E}_{n,m} >0$ because
 $\lambda_{n,m}(p)> 0$ for all $p \in \re$. For $n=1$, this gives the desired result.
\end{proof}

\begin{remark} \label{multi}
According to \eqref{UI} and \eqref{kr34x} the spectrum of the operator ${\bf H}_{m}$ consists of the ``branches" $[  {\mathcal E}_{n,m},\infty)$ where the points ${\mathcal E}_{n,m}$ are called thresholds. In cases (iii) and (iv)
\bel{kr34yy}
   {\mathcal E}_{n,m} <  {\mathcal E}_{n+1,m}
    \ee
    for all $n\in{\Bbb N}$.   Indeed, in case iii) \eqref{kr34yy}    is a consequence of the estimate  $\lambda_{n,m}(p)< \lambda_{n+1,m}(p)$ valid for all $p\in\re$ and of formula \eqref{kr34xy}. In case (iv) one has to take additionally into account that the limit of $\lambda_{n,m}(p)$ as $p\to-\infty$ is strictly smaller than that of $\lambda_{n+1,m}(p)$.
 Inequality \eqref{kr34yy}  means that there are infinitely many distinct thresholds in each of the spectra $\sigma({\bf H}_m)$, $m \in
{\mathbb Z}$, and hence in $\sigma({\bf H})$.
\end{remark}

\begin{remark} \label{multi1}
In case( iii) the multiplicity of the spectrum of all operators $\Lambda_{n,m}$ equals at least to $2$ whereas in cases (ii) and (iv) it might be equal to $1$.
\end{remark}

\section{Group velocities}
\setcounter{equation}{0}


{\bf 4.1.} 
In this subsection we obtain a
formula for the derivative $\lambda'_{n,m}(p)$, $n \in {\mathbb
N}$, $m \in {\mathbb Z}$, which yields sufficient conditions for
the monotonicity of $\lambda_{n,m}(p)$ as a function of $p$.  
Recall that the  operators
$H_m(p)$, $m \in {\mathbb N}$, $p \in \re$, were defined in the space ${\cal H}$ by formula \eqref{u1a}.

The proof of Theorem~\ref{t41} relies on integration by parts.  To prove that non-integral terms disappear at  $r=0$, we use   standard bounds on $\psi_{n,m}(r;p)$. Unfortunately, we were unable to find necessary results in the literature and therefore give their  brief proofs.

Let us consider the differential equation of Bessel type
\bel{Be1}
-  r^{-1} ( r y ')' + m^2 r^{-2} y + q(r)y=0, \quad m=0, 1,2,\ldots,
    \ee
     in a neighborhood $(0,r_{0})$ of the point $r=0$.
    If $q(r)=0$, then it has the regular $y_{0}^{(reg)}(r)=r^m$ and singular $y_{0}^{(sing)}(r)=r^{-m}$ solutions  for $m\neq 0$ and $y_{0}^{(reg)}(r)=1 $ and   $y_{0}^{(sing)}(r)=\ln r$  for $m = 0$.


  \begin{lemma} \label{l40}
Let $m\neq 0$,  and let the function $r q(r)$ belong to the class $L^{1}(0,r_{0})$. Then equation \eqref{Be1} has a solution $y^{(reg)}(r)$ satisfying the relation
\bel{Be2}
y^{(reg)}(r)=r^m + o(r^m), \quad r\to 0.
    \ee
    For its derivative, we have the bound
    \bel{Be3}
dy^{(reg)}(r)/dr=O(r^{m-1}).
    \ee
    
    Let $m=0$.  Suppose that the function $r \ln r q(r)$ belongs to the class $L^{1}(0,r_{0})$. Then equation \eqref{Be1} has a solution $y^{(reg)}(r)$ satisfying   relation
\eqref{Be2} where $m=0$.
    For its derivative, we have the bound
    \bel{Be4}
dy^{(reg)}(r)/dr=O\left(\int_{0}^r |q(s)| ds\right).
    \ee
    Moreover, if the function $r \ln^2 r q(r)$ belongs to the class $L^{1}(0,r_{0})$, then equation \eqref{Be1} has a solution $y^{(sing)}(r)$ satisfying the relation
    \bel{Be5}
y^{(sing)}(r)= \ln r  + o(1), \quad r\to 0.
    \ee
In this case any bounded solution of equation \eqref{Be1}  coincides $($up to a constant factor$)$ with the regular solution $y^{(reg)}(r)$.
 \end{lemma}
 
 \begin{proof}
 We construct the function $y^{(reg)}(r)$ as the solution of the Volterra integral equation
 \bel{Int1}
 y^{(reg)}(r)=y^{(reg)}_{0}(r)+\varkappa_{m}\int_{0}^r s (y^{(reg)}_{0}(r)y^{(sing)}_{0}(s)-y^{(reg)}_{0}(s)y^{(sing)}_{0}(r))q(s)  y^{(reg)}(s) ds
 \ee
 where $\varkappa_{m}=(2m)^{-1}$ for $m\neq 0$ and $\varkappa_{0}= -1$.
 Differentiating it twicely, we see that $y^{(reg)}(r)$ satisfies
 equation \eqref{Be1}. Equation
 \eqref{Int1} can be solved by iterations, that is
  \bel{it}
 y^{(reg)}(r)=\sum_{n=0}^\infty y^{(reg)}_{n}(r).
 \ee
  Hereby the $n^{th}$-iteration obeys the bound
 \[
| y_{n}^{(reg)}(r)|\leq \frac{C^n}{n!}r^m\Big(\int_{0}^r s  |q(s) | ds\Big)^n
 \]
 if $m\neq 0$; if $m = 0$, then $s  |q(s) |$ should be replaced by $s |\ln s| |q(s) |$. This ensures the convergence of   series \eqref{it} as well as relation  \eqref{Be2}.
 Differentiating equation \eqref{Int1} and using \eqref{Be2}, we get bounds \eqref{Be3} and \eqref{Be4} on the derivative of $y^{(reg)}(r)$.
 
 If $m=0$, we can construct the function $y^{(sing)}(r)$ as the solution of   equation
  \eqref{Int1} where the first term, $y_{0}^{(reg)}(r)$, in the right-hand side is replaced by
  $ y_{0}^{(sing)}(r)$, that is
 \[
 y^{(sing)}(r)=\ln r+\int_{0}^r s \ln (r/s) q(s)  y^{(sing)}(s) ds.
 \]
This equation can again be solved by iterations which, in particular, implies estimate \eqref{Be5}.
  \end{proof}
  
  This result can be supplemented by the following
  
   \begin{lemma} \label{l41}
Let $m\neq 0$,  and let the function $ r q^2 (r)$ belong to the class $L^{1}(0,r_{0})$.
Assume additionally that $q =\bar{q }$.  If $\psi$ is a solution of   equation \eqref{Be1} from the class $L^{2}((0,r_{0}); rdr)$, then it coincides $($up to a constant factor$)$ with the regular solution $y^{(reg)}(r)$ and hence satisfies estimates \eqref{Be2} and \eqref{Be3}.
 \end{lemma}
 
 \begin{proof}
 Let us extend the function $q(r)$     to $(r_{0},\infty)$   by zero, and let us
  consider the differential operator
$$
hy= -  r^{-1} ( r y ')' + m^2 r^{-2} y + q(r)y
$$
in the space $L^{2}({\Bbb R}_{+}; rdr)$ on domain $C_{0}^\infty({\Bbb R}_{+} )$.  
 If $q=0$, we denote this operator by $h_{0}$. 
The operator $h_{0}$ is essentially self-adjoint. To prove the same for $h$, it suffices to check that
\bel{ess}
\int_{{\Bbb R}_{+}}q^2(r) |f(r)|^2 r dr\leq \varepsilon \| h_{0}f\|^2
 +C \| f\|^2 , \q f\in C_{0}^\infty({\Bbb R}_{+} ), \q \varepsilon<1.
\ee
 Let us use the estimate
\begin{eqnarray*}
\int_{|x|\leq r_{0}}q^2(|x|) |u(x)|^2   dx\leq
\int_{|x| \leq r_{0}}q^2(|x|)     dx \, \max_{x\in {\Bbb R}^2}  |u(x)|^2 
\nonumber\\
\leq \varepsilon \int_{{\Bbb R}^2}  |(\Delta u)(x)|^2   dx 
 +C\varepsilon^{-1} \int_{{\Bbb R}^2}  |  u(x)|^2   dx, \q \forall \varepsilon>0 .
 \end{eqnarray*}
 Restricting it on the subspace of functions $u(x)=f(r)e^{i m\theta}$, we obtain estimate \eqref{ess} which implies that $h$ is essentially self-adjoint as well as $h_{0}$. Thus, equation \eqref{Be1} has at most one solution from $L^{2}((0,r_{0}); rdr)$ which is necessarily proportional to $y^{(reg)}(r)$.
  \end{proof}

Now we are in a  position to obtain a
formula for the derivative $\lambda'_{n,m}(p)$. In addition to our usual assumptions that $b(r)$ is not too singular at $r=0$, an integration-by-parts marchinery requires that $b(r)$ does not vanish too rapidly as $r\to 0$. The precise conditions are formulated rather differently    in  the cases    $m\neq 0$   and $m=0$.   We start with the first case.

\begin{theorem} \label{t41}
Let $m\neq 0$.
Suppose   that $b \in
C^3(\re_+)$  and $b(r)  
> 0$, $r \in \re_+$.
  Assume  \eqref{j1} and that $b(r)=O(e^{cr})$ for some $c>0$ as $r\to\infty$. At $r=0$
  we suppose that $b(r)=O(r^{-\gamma})$ where    $\gamma<3/2$. Moreover,
  we assume that for some $\beta< 2|m|- 1 $
  \bel{beta}
| (b(r)^{-1})^{(k)}| \leq C r^{-\beta-k}, \q k=0,1,2,3, \q r\to 0.
\ee
Put
\[
v(r)=r (r^{-1}(rb(r)^{-1})')'.
\]
  Then 
     \begin{eqnarray} 
    \lambda_{n,m}'(p) = - 2\int_0^{\infty}  r b^{-2}(r)b'(r) \psi_{n,m}'(r;p)^2 dr
    \nonumber\\
  -2^{-1}\int_0^{\infty} v' (r)  \psi_{n,m}^2 (r; p)   dr
   +2 m^2 \int_0^{\infty}   r^{-2}  b^{-1}(r)   \psi_{n,m}^2(r; p)   dr ,
   \label{sa8}\end{eqnarray}
 where the eigenfunctions $\psi_{n,m}(r;p)$ of the operator $H_{m}(p)$ are real and normalized, that is
$\|\psi_{n,m}\|  = 1$.
\end{theorem}

\begin{proof}
In view of the equation 
\bel{Bess}
 (a(r) + p)^2
\psi_{n,m} =r^{-1} ( r\psi_{n,m}' )' -m^2 r^{-2} \psi_{n,m} + \lambda_{n,m}\psi_{n,m}
\ee
we can apply to the function $\psi_{n,m}  $ the results of Lemmas~\ref{l40} and \ref{l41}
where $q(r)=(a(r) + p)^2 -\lambda_{n,m}$. Thus, Lemma~\ref{l41} implies that $\psi_{n,m}(r;p)= O(r^{|m|} )$ and 
  $\psi_{n,m}'(r;p)= O(r^{|m|-1})$    as $r\to 0$ which ensures that 
   non-integral terms     disappear at $r=0$.

To prove the same for non-integral terms corresponding to $r\to \infty$, we use   super-exponential  decay of eigenfunctions $\psi_{n,m}(r;p )$ of the  operators
$H_m(p)$. This result is valid \cite{Sh} (see also \cite{Gl}) for all one-dimensional Schr\"odinger operators with discrete spectra. In view of the condition $a(r)=O(e^{cr})$, it follows from equation
\eqref{Bess} that the derivatives $\psi_{n,m}'(r ;p)$ also decay super-exponentially.

 Let us proceed from the formula of the first order perturbation theory (known as the Feynman-Hellman formula)
\bel{F-H}
 \lambda_{n,m}'(p) = \int_0^{\infty} \frac{\partial (a(r) +
 p)^2}{\partial p} \psi_{n,m}^2 (r;p) rdr  = \int_0^{\infty}   \frac{\partial (a(r) +
 p)^2}{\partial r} \psi_{n,m}^2(r;p ) \tau(r) dr
 \ee
   where $\tau(r)=r b(r)^{-1}$.
Using that $a(r)=O(r^{1-\gamma})$  and $\tau (r)=O(r^{1-\beta})$,   we  integrate by parts and get
   \[
     \lambda_{n,m}'(p)  = - \int_0^{\infty}    (a(r) +
 p)^2   \psi_{n,m}(r;p ) ( \tau'(r)  \psi_{n,m} (r;p)  + 2  \tau(r)  \psi_{n,m}'(r;p)   )dr.
    \]
    Now it follows from   equation \eqref{Bess}
that
  \begin{eqnarray}
    \lambda_{n,m}'(p)
    = -\lambda_{n,m}(p) \int_0^{\infty}      ( \tau(r)  \psi_{n,m}^2 (r;p)   )' dr 
    +m^2 \int_0^{\infty}  r^{-2}    ( \tau(r)  \psi_{n,m}^2 (r;p)   )'dr 
     \nonumber\\
   - \int_0^{\infty}  r^{-1} ( r\psi_{n,m}' (r;p))'  ( \tau'(r)  \psi_{n,m} (r;p)  + 2\tau (r)   \psi_{n,m}'(r;p)   )dr. 
\label{sa6}\end{eqnarray}
By the  condition $ \tau (r)  \psi_{n,m}^2 (r;p) \to 0$ as $r\to 0$, the first term in the right-hand side     equals zero. In the second term we integrate by parts which yields
\[
\int_0^{\infty}  r^{-2}    ( \tau (r)  \psi_{n,m}^2 (r;p)   )'dr = 2 \int_0^{\infty}  r^{-3}   \tau (r)  \psi_{n,m}^2 (r;p)   dr
\]
because $r^{-2} \tau (r)  \psi_{n,m}^2 (r;p) \to 0$.

In the last integral in the right-hand side of \eqref{sa6}, we also integrate by parts
using that   $\psi_{n,m}'(r;p) \psi_{n,m}(r ;p) \tau'(r)\to 0$ as $r\to 0$. Thus, we have that
\begin{eqnarray}
  - \int_0^{\infty}  r^{-1} ( r\psi_{n,m}'(r;p)  )'    \tau'(r)  \psi_{n,m} (r;p)   dr
  &=& \int_0^{\infty}     \tau'(r)  \psi_{n,m}' (r;p)^2   dr 
  \nonumber\\
  &+ &   \int_0^{\infty}     v(r)  \psi_{n,m}' (r;p)   \psi_{n,m}  (r;p) dr  . 
\label{sa5}\end{eqnarray}
The last integral in the right-hand side equals
\[
- 2^{-1} \int_0^{\infty}  v'(r)  \psi_{n,m}^2 (r;p)  dr
\]
because   $  v(r) \psi_{n,m}^2(r;p) \to 0$ as $r\to 0$.
Similarly, we get that
\begin{eqnarray*}
 - 2 \int_0^{\infty}   r^{-1} ( r\psi_{n,m}'(r;p)  )'      \tau(r)   \psi_{n,m}'(r;p)   dr
 &=&- \int_0^{\infty}  r^{-2} \tau (r)        d ( r \psi_{n,m}'(r;p)^2 ) 
\nonumber\\
 &=&
  \int_0^{\infty}  r^2  (r^{-2} \tau (r))'       \psi_{n,m}'(r;p)^2   dr 
 \end{eqnarray*}
  since $\tau (r)        \psi_{n,m}'(r;p)^2\to 0$ as $r\to 0$.
 Putting the results obtained together, we arrive at representation \eqref{sa8}.
  \end{proof}
  
  \begin{follow} \label{t41cc}
 If $b'(r)\leq 0$ and $  r^{ 2} b(r) v'(r)\leq 4 m^2$ for all $r\geq 0$, then $ \lambda_{n,m}'(p)\geq 0$ for all $p\in {\Bbb R}$ and $ n$. If, moreover, one of these inequalities is strict on some interval, then $ \lambda_{n,m}'(p) >0$.
 \end{follow}

\begin{follow} \label{t41c}
 If $b(r)=b_{0} r^{-\delta}$,   $\delta\in [0, 1]$, then $\tau(r)= b_{0} ^{-1} r^{1+\delta}$,
 $v (r)= b_{0} ^{-1} (\delta^2-1) r^{\delta -1}$ and
   \begin{eqnarray*} 
    \lambda_{n,m}'(p) = 2 b_{0}^{-1} \delta \int_0^{\infty}  r^\delta  \psi_{n,m}'(r;p)^2 dr
    \\ 
   + b_{0}^{-1} (2 m^2- 2^{-1} (1-\delta)^2 (1+\delta))\int_0^{\infty}     r^{-2+\delta} \psi_{n,m}^2(r;p)   dr.
   \end{eqnarray*}
  For $b_{0}> 0$, this expression 
    is strictly positive $($so that  the functions $ \lambda_{n,m}(p)$ are strictly increasing for all $p\in {\Bbb R})$ for $m\neq 0$ since $  (1-\delta)^2 (1+\delta) \leq 1$.
    Moreover, for $\delta=1$ this result is true    for all $m\in {\Bbb Z}$.
\end{follow}

In the case $m=0$ we consider for simplicity only fields \eqref{bb}. 

\begin{pr} \label{t41cx}
 If $b(r)=b_{0} r^{-\delta}$,   $\delta\in [0, 1]$, then 
   \begin{eqnarray*} 
    \lambda_{n,0}'(p) = 2 b_{0}^{-1} \delta \int_0^{\infty}  r^\delta  \psi_{n,0}'(r;p)^2 dr
    \\ 
  - b_{0}^{-1}  2^{-1} (1-\delta)^2 (1+\delta)
  \int_0^{\infty}     r^{-2+\delta} (\psi_{n,0}^2(r;p) -\psi_{n,0}^2(0;p))  dr.
   \end{eqnarray*}
  If $b_{0}> 0$ and $\delta=1$, then  $ \lambda_{n,0}'(p)>0$
   for all $p\in {\Bbb R}$.
\end{pr}

\begin{proof}
Let us proceed again from formula \eqref{F-H}.
 We use now  that the function $\psi_{n,0}(|x|;p) $ of $x\in {\Bbb R}^2$ belongs to the Sobolev class $\mathsf{H}^2_{loc}({\Bbb R}^2)$, and therefore $\psi_{n,0}(r;p)$ has a finite limit as $r\to 0$. Thus, by Lemma~\ref{l40}  $\psi_{n,0}'(r;p)= O(r^{1-\varepsilon} )$ for any $\varepsilon>0$   as $r\to 0$. These results allow us to intergrate by parts as in the case $m\neq 0$. The only difference is with the second integral in the right-hand side of \eqref{sa5}. Now
  $v (r)= b_{0} ^{-1} (\delta^2-1) r^{\delta -1}$ and this integral  equals
  \begin{eqnarray*}
  \int_0^{\infty}     v (r)  \psi_{n,0}' (r;p)   \psi_{n,0}  (r;p) dr
 & = & 2^{-1}   \int_0^{\infty}     v(r)    d  (\psi_{n,0}^2  (r;p) - \psi_{n,0}^2  (0;p))
 \\
& = &  -  2^{-1}   \int_0^{\infty}     v'(r)      (\psi_{n,0}^2  (r;p) - \psi_{n,0}^2  (0;p))dr
 \end{eqnarray*}
because $v(r)(\psi_{n,0}^2  (r;p) - \psi_{n,0}^2  (0;p))$ as $r\to 0$.
  \end{proof}

\medskip

{\bf 4.2.}
 In this subsection we  show that for linear potentials, that is for magnetic fields not depending on $r$, all eigenvalues $\lambda_{n,0}(p)$, $n \in {\mathbb
N}$, of the operator $H_0(p)$ are not monotonous functions of
$p\in\re$. We follow closely the proof of the first part of Proposition~\ref{t33}. However   we now use that eigenfunctions of the harmonic oscillator decay faster than any power of $r^{-1}$ at infinity (actually, they decay super-exponentially).

\begin{pr} \label{p41} Assume that  for sufficiently large $r$
  \bel{kr31}
    b (r) =  b_{0}>0.
    \ee
  Then, for all $n \in {\mathbb N}$, some $\gamma_{n}>0$ and
sufficiently large $k>0$,  we have
    \bel{san30}
    \lambda_{n,0}(-k) \leq ( 2n-1)b_{0} -\gamma_{n}k^{-2}.
    \ee
\end{pr}

\begin{proof}
Let $\zeta$ be the same function as in the proof of the first part of Proposition~\ref{t33}. We set $\rho_{k}=b_{0}^{-1} k$, $\gamma_{k}=2 b_{0}^{1/2} k^{-1}$ and define the functions $\zeta(r;k)$ and $\varphi_{j}(r;k)$ by formulas \eqref{j9b} and \eqref{j9}, respectively. It suffices to check that
 \bel{san35}
   \langle L_0(-k)\varphi(k),\varphi (k) \rangle
    \leq 2n-1   -\gamma_{n}k^{-2}.
    \ee
    for sufficiently large $k$ and  all normalized functions from subspace \eqref{j9a}.
    Let us proceed from formula \eqref{j5x1}.
Since the functions $\psi_{j}(x;k)$  decay faster than any power of $|x|^{-1}$ as  $|x| \to\infty$,   the term $o(1)$ in  \eqref{j19a} is actually $O(k^{-\infty})$.   Similarly, estimate \eqref{j12} can be formulated in a more precise form as
\bel{psi1}
\|\psi(k) \zeta'(k)\|^2= O(k^{-\infty}).
\ee
Since $r\leq 2^{-1} 3k$ on the support of $\varphi(k)$, we have that 
\bel{psi2}
\| r^{-1} \psi(k)  \|^2 \geq (2/3)^2 k^{-2}.
\ee
Now   function \eqref{psi} is zero if $r$ and $k$ are large enough. Therefore equation \eqref{j14} yields the exact equality
 \begin{eqnarray*} 
 {\rm Re} \,  \langle  -\psi''(k) + (b_{0}r -k)^2 \psi (k), \psi (k)\zeta^2(k)\rangle   =
b_{0} \sum_{j,l=1}^n (2j-1) c_{j}\bar{c}_{l}  \langle  \varphi _{j}(k), \varphi_{l} (k)\rangle   
    \end{eqnarray*}
(cf. \eqref{j5x2}). Up to terms $O(k^{-\infty})$, the right-hand side here is estimated by $b_{0}(2n-1)$. Together with \eqref{psi1} and \eqref{psi2}, this implies estimate \eqref{san35}.
  \end{proof}

    Combining relations \eqref{j4} and \eqref{san30}, we see that the eigenvalues
    $\lambda_{n,0}(p)$ tend as $p\to -\infty$ to their limits $(2n-1) b_{0}$ from below. On the other hand, according to \eqref{kr4} $\lambda_{n,0}(p)\to \infty$ as $p\to \infty$. Thus, all functions  $\lambda_{n,0}(p)$ have necessarily local minima. We can obtain an additional information using the following elementary
    
    \begin{lemma} \label{f41h}
Suppose that \eqref{kr31} is satisfied for all $r>0$ and that $a(r)=b_{0} r$. Then
 \bel{harm}
\lambda_{n,m}(0)= 2 b_{0} ( 2n-1+ |m|)
\ee
  for all $n \in {\mathbb N}$ and $m \in {\mathbb Z}$.  
\end{lemma}
    
 \begin{proof}
 Let us consider the two-dimensional harmonic oscillator ${\bf T}=-\Delta + b^2_{0} ( x_{1}^2+x_{2}^2)$.
 Separating the variables $x_{1}$, $x_{2}$, we see that its spectrum consists of the eigenvalues
 $2 b_{0} (l_{1}+l_{2}-1)$ where $l_{1}, l_{2} \in {\mathbb N}$. It follows that the operator ${\bf T}$ has the eigenvalues $2b_{0}j$, $j \in {\mathbb N}$, of multiplicity $j$. On the other hand, separating the variables in the polar coordinates, we see that the spectrum of ${\bf T}$ consists of the eigenvalues $\lambda_{n,m}(0)$ of the operators $H_{m}(0)$. For the proof of \eqref{harm} we
  take into account that all eigenvalues $\lambda_{n,m}(0)$ are simple and that
 $\lambda_{n,m+1}(0) > \lambda_{n,m}(0)$ for all $n$ and $m\geq 0$. Clearly, the operator $H_{0}(0)$ has an eigenvalue $2b_{0}j$ if and only if its multiplicity $j$ is odd. This gives formula \eqref{harm} for $m=0$.
  We shall show that for every $j\in {\mathbb N}$
 \bel{harm1}
\lambda_{1,j-1}(0)=\lambda_{2,j-3}(0)=\ldots=\lambda_{2,-j+3}(0)=\lambda_{1,-j+1}(0)=2b_{0}j
\ee
which is equivalent to  formula \eqref{harm} for all $m$.  
  Let us choose some $j_{0}$ and suppose that \eqref{harm1} holds for all  $j\leq j_{0}$. Then we check it for $j=j_{0}+1$. First we remark that if an  operator  $H_{m}(0)$ for some $m>0$ has $n$ eigenvalues in the interval $[2b_{0}, 2b_{0}(j_{0}+1)]$, then the operator  $H_{m-1}(0)$ has at least $n$ eigenvalues in the interval $[2b_{0}, 2b_{0} j_{0}]$. Then using \eqref{harm1} for $j\leq j_{0}$, we see that if an  operator  $H_{m}(0)$ has
   the eigenvalue $2b_{0}(j_{0}+1)$, then necessarily the operator  $H_{m-1}(0)$   has   the eigenvalue $2b_{0}j_{0} $. Therefore according to \eqref{harm1} for $j=j_{0}$, only the operators $H_{m}(0)$  with $m=j_{0}, j_{0}-2,\ldots, -j_{0}+2, -j_{0}$ might have  the eigenvalue $2b_{0}(j_{0}+1)$. There are $j_{0}+1$ of such operators and the multiplicity of this eigenvalue equals $j_{0}+1$. Thus, all the operators $H_{m}(0)$  for $m=j_{0}, j_{0}-2,\ldots, -j_{0}+2, -j_{0}$ and only for such $m$ have  the eigenvalue $2b_{0}(j_{0}+1)$. This proves \eqref{harm1}   for   $j_{0}+1$.
  \end{proof}

   Comparing this result   with \eqref{j6}, we see that, for potentials $a(r)=b_{0} r$,
   \[
   \lim_{p\to-\infty} \lambda_{n,m}(p)=b_{0}(2n-1)< 2 b_{0}( 2n-1+ |m|) = \lambda_{n,m}(0).
   \]
   Together with \eqref{san30}, this implies that  the functions $\lambda_{n,0}(p)$ have   negative local minima.

     Thus, we get the following

\begin{theorem} \label{f41}
Under the hypotheses of Proposition~$\ref{p41}$ the eigenvalues
$\lambda_{n,0}(p)$, $n \in {\mathbb N}$, of the operator $H_0(p)$
are not monotonous functions of $p \in \re$. Moreover, if \eqref{kr31} is satisfied for all $r>0$ and $a(r)=b_{0} r$, then the functions $\lambda_{n,0}(p)$ lose their monotonicity for $p<0$.
\end{theorem}

We do not know how many minima have the functions $\lambda_{n,0}(p)$.

The problem of monotonicity of  the eigenvalues
$\lambda_{n,0}(p)$ for fields $b(r)=b_{0}r^{-\delta}$ where $\delta\in (0,1)$ remains also open.

\medskip

{\bf 4.3.} 
In a somewhat similar situation the break down of monotonicity of group velocities was exhibited in \cite{H}. In this paper one considers the Schr\"odinger operator ${\bf H}^{(N)} =
-\frac{\partial^2}{\partial x^2}+ \left(i\frac{\partial}{\partial
y} - bx\right)^2$ with constant magnetic field $b>0$, defined on
the
semi-plane $\left\{(x,y) \in \rd \, : \, x>0\right\}$ with the Neumann boundary condition at $x=0$.
Let $H^{(N)}(p)=-d^2/dx^2 + (bx+p)^2$, $p \in \re$, be the self-adjoint operator
in the space $L^2(\re_+)$ corresponding to the   boundary condition $u'(0)=0$.
 Then the operator ${\bf H}^{(N)}$ is unitarily
equivalent under  the partial
Fourier transform with respect to $y$, to the direct integral
$\int_{\re}^{\oplus} H^{(N)}(p) dp$. It is shown in \cite[Section 4.3]{H} that
 the lowest eigenvalue  $\mu_1(p)$  of $H^{(N)}(p)$  is not monotonous for $p<0$. This follows from the inequality $ \mu_1'(0) > 0$ proven\footnote{Note
that in \cite{DH} and \cite{H} the parameter $p$ is chosen with
the opposite sign.} in \cite{DH} and the relations
    \bel{san41}
      \lim_{p \to -\infty} \mu_1(p)= \mu_1(0)  = b.
    \ee

Our proof of non-monotonicity of the functions  $\lambda_{n,0}(p)$ is essentially different since in contrast with \eqref{san41}  we have $\lim_{p \to -\infty} \lambda_{n,0}(p) <
\lambda_{n,0}(0)$. 

\section{Asymptotic time evolution }
\setcounter{equation}{0}


{\bf 5.1.} 
Combined with the stationary phase method, the spectral analysis of the operators
${\bf H} ={\bf H} (a)$     allows us to find the asymptotics for large $t$
of solutions $u(t)=\exp(-i{\bf H}  t)u_0$ of the time dependent Schr\"odinger equation. It follows
from (\ref{MH}) that
\[
  \overline{\exp(-i{\bf H} (a) t) u_0}=\exp(i{\bf H}  (-a) t) \overline{  u_0}.
\]
Therefore it suffices to consider the case $a(r)\to+\infty$. 
Moreover, on every subspace ${\goth H}_m$ with a fixed
magnetic quantum number $m$, the problem reduces to the asymptotics of the function 
 $u(t)=\exp(-i {\bf H}_m  t)u_0$.
 
 Let us proceed from decomposition (\ref{UI}). Suppose that ${\cal F}_{m} u_{0}\in\ran \Psi_{n,m} $. Then (see (\ref{Psi}))   
 \begin{equation}\label{eq:InV}
( {\cal F}_{m} u_{0} ) (r,p)=\psi_{n,m} (r,p)f(p) 
\end{equation}
where $f= \Psi_{n,m}^* {\cal F}_{m} u_{0}$ and $u(t)= {\cal F}_{m}^*\Psi_{n,m} e^{-i\Lambda _{n,m}t} f$, that  is
\begin{equation}\label{eq:TE}  
     u_{n,m}(r,x_{3},t)=(2\pi)^{-1/2}\int_{-\infty}^\infty e^{ip x_{3} -i\lambda_{n,m}(p) t}\psi_{n,m}(r,p)f(p)dp.
\end{equation}

The analytic function $\lambda_{n,m}''(p)$ might have only a countable set of zeros $p_{n,m,l}$ with possible accumulations at $\pm\infty$ only. The function $\lambda_{n,m}'(p)$ is monotone on every interval $ (p_{n,m,l}, p_{n,m,l+1})$ and takes there all values between $\lambda_{n,m}'(p_{n,m,l})=:\alpha_{n,m,l}$ and $\lambda_{n,m}'(p_{n,m,l+1})=:\beta_{n,m,l}$.  
 We consider the asymptotics of integral  \eqref{eq:TE}  on each of the subspaces $L^2 (p_{n,m,l}, p_{n,m,l+1})$ separately.
 Let us set $\gamma=x_{3} t^{-1}$.
  First we suppose that $f\in C_{0}^\infty (p_{n,m,l}, p_{n,m,l+1})$. 
The stationary points of   integral \eqref{eq:TE}  are determined by the equation
\begin{equation}\label{eq:stp}  
      \lambda_{n,m}^\prime(p)=\gamma.
\end{equation}

If $\gamma\not\in (\alpha_{n,m,l}, \beta_{n,m,l})$, it does not have  solutions from the interval $(p_{n,m,l}, p_{n,m,l+1})$. Therefore integrating directly by parts, we find that function \eqref{eq:TE} decays in this region of $x_{3}/t$ faster than any power of $(|x_{3}|+|t|)^{-1}$ (and $r$). 
If $\gamma\in (\alpha_{n,m,l}, \beta_{n,m,l})$,  then on the interval $(p_{n,m,l}, p_{n,m,l+1})$ equation \eqref{eq:stp} has a unique solution which we denote by $\nu_{n,m,l} (\gamma)$. Let us   set
\[
      \Phi_{n,m,l}(\gamma)=\nu_{n,m,l} (\gamma)\gamma -\lambda_{n,m}(\nu_{n,m,l}(\gamma))
\]
and denote by $\chi_{n,m,l}$  the characteristic function of the interval $  (\alpha_{n,m,l}, \beta_{n,m,l})$.
  For   $\gamma$ from this interval, we apply the stationary phase method to integral \eqref{eq:TE} which yields
\begin{eqnarray}\label{eq:TE1}  
     u (r,x_{3},t)= \tau_{n,m,l}^{(\pm)} e^{i\Phi_{n,m,l}(\gamma)t }
\psi_{n,m} (r,\nu_{n,m,l}(\gamma)) |\lambda_{n,m}'' (\nu_{n,m,l} (\gamma)) |^{-1/2}
\nonumber\\
\times f(\nu_{n,m,l}(\gamma)) \chi_{n,m,l} (\gamma) |t|^{-1/2}
+ u_\infty(r,x_{3},t), \q \gamma=x_{3} t^{-1}, \q t\to\pm \infty,
\end{eqnarray}
where 
$\tau_{n,m,l}^{(\pm)}= e^{ \mp \pi i \sgn (\lambda''_{n,m}(p))/4}$ for $p\in (p_{n,m,l}, p_{n,m,l+1})$ and
\begin{equation}\label{eq:TE1r} 
\lim_{t\rightarrow\pm\infty }
\| u_\infty(\cdot,t)\|= 0.
\end{equation}
 Since the norm in   the space  ${\goth H}$
of the first term in the right-hand side of (\ref{eq:TE1}) equals the norm of $f$ in the space $ L^2 (p_{n,m,l}, p_{n,m,l+1})$,
 asymptotics (\ref{eq:TE1}) extends   to all functions (\ref{eq:InV}) with an arbitary 
 $f \in L^2 (p_{n,m,l}, p_{n,m,l+1})$.
 Thus, we have proven

\begin{theorem}\label{time}
Assume \eqref{kr8} and \eqref{j1}. 
Let  $u(t)=\exp(-i {\bf H}_m t) u_0$ where $u_0$ satisfies $(\ref{eq:InV})$
with $f\in L^2 (p_{n,m,l}, p_{n,m,l+1})$. Then the
asymptotics as $t\rightarrow \pm \infty$ of this function  
 is given by relations $(\ref{eq:TE1})$,
$(\ref{eq:TE1r})$.
\end{theorem}

Of course asymptotics $(\ref{eq:TE1})$, $(\ref{eq:TE1r})$ extends automatically to all 
$f \in L^2({\Bbb R})$ with compact support and to linear of functions
$\psi_{n,m} (r,p)f_{n}(p)$ over different $n$.

By virtue of formulas $(\ref{eq:TE1})$, $(\ref{eq:TE1r})$ a quantum particle in   magnetic field \eqref{magnfi} remains localized in the $(x_{1},x_{2})$-plane  but propagates in the $x_{3}$-direction.   If $f\in L^2 (p_{n,m,l}, p_{n,m,l+1})$, then a particle ``lives" as $|t|\to\infty$ in the region where $x_{3}\in (\alpha_{n,m,l} t,
\beta_{n,m,l} t)$. In particular, if $  \lambda' (p)>0$ ($  \lambda' (p)<0$) for $p\in (p_{n,m,l}, p_{n,m,l+1})$, then a particle propagates in the positive (negative) direction as $t\to+\infty$. Thus, according to Corollary~\ref{t41c} if $b(r)=b_{0}r^{-\delta}$, $\delta\in [0,1]$, $b_{0}>0$, then a particle with the magnetic quantum number $m\neq 0$ propagates always in the positive direction of the $x_{3}$-axis. If $\delta =1$, then this   result remains true from all $m$. On the contrary, if $\delta =0$ and $m=0$, then a particle will   propagate in a negative direction for some interval of momenta $p$.

\medskip

{\bf 5.2.}
 Theorem~\ref{time} implies the existence of asymptotic velocity in the $x_{3}$-direction.
 The corresponding operator is defined by the equation (cf. \eqref{UI})
\[
   {\bf H}_{m}'      =\bigoplus_{n\in {\Bbb N}} {\cal F}_{m} ^* \Psi _{n,m} \Lambda_{n,m}' \Psi_{n,m}^*  {\cal F}_{m} , 
\]
where $\Lambda_{n,m}' $ are the operators of multiplication by the functions $\lambda_{n,m}' (p)$.
To put it differently, the operator ${\bf H}_{m}'  $ acts as  multiplication by   $\lambda_{n,m}' (p)$ in the spectral representation of the operator ${\bf H}_{m}  $ where it   acts as  multiplication by  the functions $\lambda_{n,m}(p)$.

\begin{pr}\label{asve}
Assume \eqref{kr8} and \eqref{j1}. 
Then, for an arbitrary bounded function ${\mathcal Q}$,
\bel{UIy}
 \slim_{|t| \to \infty} \exp{(i{\bf H}_m t)} {\mathcal
Q}(x_{3}/t) \exp{(-i{\bf H}_m t)} = {\mathcal Q}\left({\bf
H}'_m\right)
\ee
$($in particular, the strong limit in the left-hand side exists$)$.
\end{pr}

\begin{proof}
We shall check that for all $u_{0}\in{\cal H}_{m}$
\bel{UIz}
 \lim_{|t| \to \infty}   \| {\mathcal
Q}(x_{3}/t) \exp{(-i{\bf H}_m t)} u_{0} - \exp{(-i{\bf H}_m t)}  {\mathcal Q} ({\bf
H}'_m) u_{0} \| =0
\ee 
which is equivalent to relation \eqref{UIy}.
 Remark that if $u_{0}$ satisfies    \eqref{eq:InV}, then
\bel{UIzz}
 ( \Psi_{n,m}   {\mathcal Q}  ({\bf H}'_m) u_0) (r,p)
 =\psi_{n,m} (r,p){\mathcal Q}  (\lambda_{n,m}' (p)) f(p).
\ee
 It suffices to prove \eqref{UIz} on a dense set of elements $u_{0}$ such that equality  \eqref{eq:InV} is true with  $f\in L^2 (p_{n,m,l}, p_{n,m,l+1})$.   Applying the operator ${\mathcal
Q}(x_{3}/t)$  to asymptotic relation \eqref{eq:TE1}, we see that the asymptotics of 
${\mathcal Q}(x_{3}/t) \exp{(-i{\bf H}_m t)} u_{0}$ is given again by formula \eqref{eq:TE1} where
the function $f(\nu_{n,m,l}(\gamma))$ in the right-hand side is replaced by the function
${\mathcal Q} (\gamma)f(\nu_{n,m,l}(\gamma))$. Similarly, it follows from Theorem~\ref{time} and relation \eqref{UIzz} that the asymptotics of  $ \exp{(-i{\bf H}_m t)}  {\mathcal Q} ({\bf
H}'_m) u_{0} $ is given by formula \eqref{eq:TE1} where
the function $f(\nu_{n,m,l}(\gamma))$ in the right-hand side is replaced by the function
${\mathcal Q} (\lambda'_{n,m} (\nu_{n,m,l}(\gamma)))f(\nu_{n,m,l}(\gamma))$. So for the proof  of \eqref{UIz},  it remains to take equation \eqref{eq:stp} into account.
\end{proof}

Relation \eqref{UIy} shows that ${\bf H}_m'$ can naturally be interpreted as the operator
of asymptotic velocity in the $x_{3}$-direction.

Similar results concerning the Iwatsuka model (see \cite{i} or
\cite{CFKS}) have been obtained in \cite{mp}.

Numerous useful discussions with Georgi Raikov as well as a financial support
 by the Chilean Science Foundation {\em Fondecyt} under
Grant  7050263 are gratefully   acknowledged.

\noindent D\'epartement de math\'ematiques, 
Universit\'e de Rennes I,
Campus Beaulieu, 35042  Rennes, FRANCE,  
yafaev@univ-rennes1.fr
\end{document}